\documentclass[12pt]{article}
\usepackage{amssymb,amsmath}
\renewcommand{\tilde}{\widetilde}
\renewcommand{\hat}{\widehat}
\renewcommand{\phi}{\varphi}
\newtheorem{theorem}{Theorem}[section]
\newtheorem{df}[theorem]{Definition}
\newtheorem{lm}[theorem]{Lemma}
\newtheorem{cor}[theorem]{Corollary}
\newtheorem{rem}[theorem]{Remark}
\newtheorem{remexs}[theorem]{Remarks and Examples}
\newtheorem{prop}[theorem]{Proposition}

\newtheorem{conj}[theorem]{Conjecture}
\newcommand{\proof}{{\noindent \rm\bf Proof. }}
\newcommand{\qed}{{$\diamond$ \bigskip}}
\newcommand{\benum}{\begin{enumerate}}
\newcommand{\ennum}{\end{enumerate}}
\newcommand{\N}{\mathbb{N}}

\newcommand{\gen}[1]{\mbox{$\langle #1 \rangle$}}
\newcommand{\gena}[1]{\mbox{$\{ #1 \}$}}
\newcommand{\spin}[1]{\mbox{$\langle #1 \rangle$}}
\newcommand{\sgn}{{\mbox{\rm sgn} \: }}
\newcommand{\la}{\lambda}
\DeclareMathOperator{\Irr}{Irr}
\DeclareMathOperator{\GL}{GL}
\DeclareMathOperator{\PSL}{PSL}
\title{\bf
Critical classes, \\
Kronecker products of spin characters, \\
and the Saxl conjecture
}
\author{Christine Bessenrodt\\
\small Institute for Algebra, Number Theory and Discrete Mathematics\\
\small Leibniz Universit\"at Hannover\\[-0.8ex]
\small Welfengarten 1, D-30167 Hannover, Germany\\
\small\tt bessen@math.uni-hannover.de}

\date{April 3, 2017}

\begin{document}

\maketitle


\begin{abstract}
Highlighting the use of critical classes,
we consider constituents in
Kronecker products, in particular of
spin characters of the double covers
of the symmetric and alternating groups.
We apply results from the spin case to find
constituents in  Kronecker products
of characters of the symmetric groups.
Via this tool,
we make progress on the Saxl conjecture;
this claims that for a triangular number $n$, the square of the
irreducible character of the symmetric group $S_n$
labelled by the staircase
contains all irreducible characters
of $S_n$ as constituents.
We find a large number of constituents
in this square which were not detected
by other methods.
Moreover, the investigation of Kronecker products
of spin characters inspires a spin variant of Saxl's conjecture.
\end{abstract}
\footnotetext{}

\section{Introduction}

Even though tensor products of complex representations and the
corresponding Kronecker products of characters of the symmetric group $S_n$ have been studied for a long time,
their decomposition is still an elusive central open problem
in the field.
A new benchmark for this is a conjecture by
Heide, Saxl, Tiep and Zalesskii \cite{HSTZ}
which says that for any $n\neq 2,4,9$ there is always an irreducible
character of $S_n$ whose square contains all irreducible characters.
For triangular numbers  an explicit candidate was suggested by Saxl
in 2012; denoting the irreducible character of $S_n$ associated to a partition $\la$ of $n$ by $[\la]$, his conjecture is the following.
\smallskip

{\bf Saxl's Conjecture } Let $\rho_k=(k,k-1,\ldots,2,1)$ be the staircase partition of $n=k(k+1)/2$. Then the Kronecker square $[\rho_k]^2$ contains all irreducible characters of $S_n$ as constituents.
\medskip

This has inspired a lot of recent research.
In particular, in the work of Pak, Panova and Vallejo \cite{PPV}
and Ikenmeyer \cite{I} many constituents of the square $[\rho_k]^2$ have been identified,
notably those to hooks and to partitions comparable
to the staircase in dominance order.

\smallskip

Here, we will take a very different approach and show that results on the spin characters of a double cover $\tilde{S}_n$ of the symmetric group $S_n$ can be fruitfully applied  towards this problem.
Indeed, it will turn out to be useful to consider also the
spin characters  of the double cover $\tilde{A}_n$
of the alternating group~$A_n$.
An important link is made as a consequence of an identity found
in an earlier investigation of homogeneous Kronecker products.
While in  \cite{BK} it was shown that there are no nontrivial homogeneous products of $S_n$-characters,
for the double covers $\tilde{S}_n$ of the symmetric groups
nontrivial homogeneous spin products do occur
for all triangular numbers~$n$ \cite{BKspin}.
In \cite{Bmix}, also nontrivial homogeneous
mixed products of complex characters
for the double covers $\tilde{S}_n$ were found,
i.e., products of a non-faithful character of  $\tilde{S}_n$
with a spin character.
This information on nontrivial homogeneous Kronecker products
is crucially used here towards obtaining many constituents
in the square $[\rho_k]^2$.

\smallskip

We give a brief overview on the following sections.
In Section~2, we collect the information and notation on the
irreducible characters for the symmetric and alternating
groups and their double covers. Here we already draw
attention to special conjugacy classes of these groups, and we prove
a general lemma that uses critical properties of conjugacy classes
for identifying constituents in certain Kronecker products.
In Section~3 we first recall properties of the basic spin character
as well as an important link between the
faithful and non-faithful characters of $\tilde{S}_n$
labelled  by the staircase partition $\rho_k$; as an
easy application we obtain all hook characters in $[\rho_k]^2$.
Section~4 starts with a short argument for the criterion given
in \cite[``Main Lemma'']{PPV},
and then moves on to the main results on
constituents in certain spin character products
(Theorems~\ref{spinMain} and~\ref{spinMainbound}).
These results are then applied towards the Saxl conjecture.
A powerful new criterion for constituents in $[\rho_k]^2$ is given
in Corollary~\ref{cor:newcriterionSaxl};
its usefulness is illustrated by providing several families of
constituents. For characters to 2-part partitions this also involves
the almost strict unimodality of the number of $k$-bounded
strict partitions.
The final Section~5 discusses a spin variant of Saxl's conjecture,
involving the ``spin staircase'' $(2k-1,2k-3,\ldots,3,1)$.
Also spin variants of the conjecture by Heide, Saxl, Tiep and
Zalesskii are presented.


\section{Preliminaries}\label{prelim}

We denote by $P(n)$ the set of partitions of $n$, i.e.,
weakly decreasing sequences of nonnegative integers summing to~$n$.
For a partition $\lambda \in P(n)$, $l(\lambda)$ denotes its length,
i.e., the number of positive parts of $\lambda$.
The set of partitions of $n$ into odd parts only is denoted by $O(n)$,
and the set of partitions of $n$ into distinct parts is denoted by
$D(n)$.
We write $D^+(n)$ resp.\ $D^-(n)$ for the sets of partitions $\lambda$
in $D(n)$ with $n - l(\lambda)$ even resp.\ odd; the partition
$\lambda$ is then also called even resp.\ odd.

We write $S_n$ for the symmetric group on $n$ letters, and
$\tilde{S}_n$ for one one of its double covers;
so $\tilde{S}_n$ is a non-split extension
of $S_n$ by a central subgroup $\gen{z}$ of order~2.
It is wellknown that the representation theory
of these double covers is 'the same' for all
representation theoretical purposes.
For $\lambda \in P(n)$, we write $[\lambda]$ for the corresponding
irreducible character of $S_n$; this is identified with the
corresponding non-faithful character of $\tilde{S}_n$.
When we evaluate $[\la]$ on an element of $S_n$ of cycle type $\mu \in P(n)$,
we simply write $[\la](\mu)$ for the corresponding value.
For background on the representation theory of the symmetric groups, the reader is referred to \cite{J,JK}.

The spin characters of $\tilde{S}_n$
are those that do not have $z$ in their kernel.
For  an introduction to the properties of spin characters
resp.\ for some results we will need in the sequel
we refer to  \cite{HH}, \cite{Mo62}, \cite{Mo65},
\cite{Schur}, \cite{Stembridge}.
Below we collect some of the necessary notation
and some results from \cite{Stembridge}
that are crucial in later sections.
For $n\le 3$, the irreducible $\tilde S_n$-characters are lifted from $S_n$.
Hence we will always assume that $n\ge 4$.

\medskip

Let $\lambda \in P(n)$.
Then the conjugacy class of elements in $\tilde S_n$
projecting to cycle type $\lambda$
splits into two $\tilde S_n$-conjugacy classes
if and only if $\lambda \in O(n) \cup D^-(n)$,
otherwise it does not split.
When it splits, a specific labelling $C^\pm_\lambda$
for the two $\tilde S_n$-classes in $C_{\lambda}$
is made; we leave out the details here (see \cite{Schur, Stembridge}), and just note that for a given spin character
the values on the two classes differ only by a sign.
Note that any spin character vanishes on the non-split classes,
thus only the values on the split classes will be considered.

Let $\sgn$ denote the sign character of~$\tilde S_n$, inflated from the sign character of $S_n$.
A character $\chi$ of~$\tilde S_n$ is called {\em self-associate} if
$\sgn \cdot \chi = \chi$, otherwise we have a pair
of {\em associate} characters,  $\chi \ne \sgn \cdot \chi$.

In 1911, Schur has proved the following classification result \cite{Schur}, giving a complete list of irreducible complex spin characters of~$\tilde S_n$.
\\
For each $\lambda \in D^+(n)$ there is a self-associate spin character
$\gen{\lambda}$, and for each $\lambda \in {\cal D}^-(n)$ there is a pair
of associate spin characters $\gen{\lambda}_+$ and $\gen{\lambda}_-$.
When we want to consider a spin character labelled by $\lambda$,
and it is not specified whether $\lambda$ is in $D^+$ or~$D^-$, we write
$\spin{\lambda}_{(\pm)}$.
\\
The spin characters labelled by $ \lambda = (\ell_1 , \dots , \ell_m)\in D(n)$
take the following values on $\sigma_{\alpha} \in C_{\alpha}^+$:
\[ \begin{array}{rcll}
\spin{\lambda}_+ (\sigma_{\alpha}) &=& \spin{\lambda}_-(\sigma_{\alpha})
      & \mbox{for } \alpha \in {O}(n), \lambda \in {D}^-(n)\\[5pt]
\spin{\lambda}_{(\pm)} (\sigma_{\alpha}) &=& 0
      & \mbox{for } \alpha \in {D}^-(n), \lambda \neq \alpha\\[5pt]
\spin{\lambda}_+ (\sigma_{\lambda}) &=& - \spin{\lambda}_-(\sigma_{\lambda})
      = i^{(n-m+1)/2} \sqrt{\frac{\prod_j \ell_j}{2}}
      & \mbox{for } \lambda  \in {D}^-(n)\\
\end{array}\]
The values $\spin{\lambda}_{(\pm)} (\sigma_{\alpha})$ for $\alpha \in {O}(n)$
are integers determined by a recursion rule akin to the Murnaghan-Nakayama rule, which is due to Morris \cite{Mo62,Mo65}.

\smallskip

Furthermore, we define
$$
\hat{\gen{\lambda}}  =  \left\{
\begin{array}{ll}
\gen{\lambda} & \mbox{if } \lambda \in D^+(n)\\
\gen{\lambda}_+ + \gen{\lambda}_-& \mbox{if } \lambda \in D^-(n)\\
\end{array}
\right.
$$

An important role in the theory is taken by the
{\em basic spin characters} $\gen{n}_{(\pm)}$.
Their values on an element $\sigma\in \tilde{S}_n$
projecting to cycle type $\alpha$ are determined explicitly
as follows.
\\
When $n$ is odd,
$$
\gen{n}(\sigma)= 2^{(l(\alpha)-1)/2}
\quad \text{ for } \alpha \in O(n)\:,$$
and when $n=2k$ is even,
$$\gen{n}_{\pm}(\sigma)=
\left\{ \begin{array}{ll}
2^{(l(\alpha)-2)/2} & \text{for } \alpha \in O(n)\\
\pm i^k \sqrt{k} & \text{for } \alpha=(n)
\end{array} \right. \:.
$$
All other values of $\gen{n}_{(\pm)}$ are zero.

\medskip

We will also need some information about the characters of the alternating groups~$A_n$ and their double covers $\tilde{A}_n$.

The classification of the irreducible $A_n$-characters is derived
from that for $S_n$. We obtain all irreducible
characters of $A_n$ as constituents in the restriction
of the characters $[\lambda]$ as follows (see \cite{JK}).
\\
Let $\mu \in P(n)$, and let
$\mu'$ be the transposed partition to $\mu$.
\\
When $\mu \neq \mu'$,
$[\mu]\downarrow_{A_n}=[\mu']\downarrow_{A_n}={\{\mu\}}=\{\mu'\} \in \Irr(A_n)$.
\\
When $\mu=\mu'$,
$[\mu]\downarrow_{A_n} = {\{\mu\}_{+}} \: + \:
{\{\mu\}_{-}}$; the characters $\{\mu\}_{\pm}$ are
conjugate irreducible characters of $A_n$.
Let $h(\mu)=(h_1,\ldots,h_k)$ be the
partition of principal hook lengths $h_1,\ldots,h_k$ in $\mu$,
i.e.\ here,  $h_j=2(\mu_j-j)+1$.
The two characters $\{\mu\}_{\pm}$ differ only on the two conjugacy classes of cycle type $h(\mu)$; note that
$[\mu]$ takes the value $e_\mu=(-1)^{(n-k)/2}$ on elements
of this cycle type.
With $\sigma^\pm_{h(\mu)}$ being appropriate
representatives of the conjugacy classes of~$A_n$ of this cycle type,
the character values are given as
$$
\{\mu\}_+(\sigma^\pm_{h(\mu)}) =
\frac 12 \left( e_\mu \pm \sqrt{e_\mu\prod_{j=1}^k h_j} \right)
$$
and similarly, with interchanged signs, for  $\{\mu\}_-$.
In particular, these two conjugate characters
are the only irreducible characters that differ on the elements $\sigma^\pm_{h(\mu)}$.
\medskip

One also obtains the classification of irreducible spin characters of
$\tilde A_n$ from that of the spin characters of $\tilde S_n$
(see \cite{HH, Schur, Stembridge}).
\\
For each $\la\in D^-(n)$, the restriction to $\tilde A_n$
gives one irreducible spin character
$$\gen{\la}_+\downarrow_{\tilde A_n} =
\gen{\la}_-\downarrow_{\tilde A_n} = \gen{\gen{\la}}\:.$$
Dually, for each $\la=(\ell_1,\ldots,\ell_m) \in D^+(n)$,
the restriction to $\tilde A_n$
gives two conjugate irreducible spin characters
$$\gen{\la}\downarrow_{\tilde A_n} =
 \gen{\gen{\la}}_+ + \gen{\gen{\la}}_-\:.$$
If $\sigma\in \tilde A_n$ projects to cycle type $\la$,
then for the difference of the values on $\sigma$ we have
(with the sign depending on the choice of associates)
$$
\Delta^\la(\sigma)= \gen{\gen{\la}}_+(\sigma) - \gen{\gen{\la}}_- (\sigma)
 = \pm i^{(n-m)/2}\sqrt{\prod_{j=1}^m \ell_j}\:;
$$
if $\sigma\in \tilde A_n$ does not project to type $\la$,
then $\Delta^\la(\sigma)=0$ \cite{Schur}.
Let $\sigma$ project to type $\la$.
If $\la \in D\cap O(n)$,   and
$\tau \in \tilde{A}_n$
is $\tilde{S}_n$-conjugate to $\sigma$
but not in $\sigma^{\tilde{A}_n}$, then
$\Delta^\la(\sigma)= -\Delta^\la(\tau)$
(see \cite[remark after 7.5]{Stembridge}).
Note that $\gen{\la}(\sigma)=0$ if $\la \in D^+(n)\setminus(D\cap O(n))$;
thus we can compute all the character values also in the case
where $\la\in D^+(n)$.
\\
For elements projecting to other types, the values of the
two conjugate characters $\gen{\gen{\la}}_\pm$ are the same;
in particular, the characters $\gen{\gen{\la}}_\pm$ vanish on classes
projecting to a cycle type $\mu\ne \la$ that is not in $O(n)$.
We emphasize that when $\la \in D\cap O(n)$,
both of these two conjugate characters
have different values on the two $\tilde{A}_n$-classes that come from
one $\tilde{S}_n$-class projecting to
cycle type $\lambda$, and they are the only irreducible spin
characters of $\tilde A_n$ with this property.
\smallskip

We also note the following special situation.
When $\la=\la'$, the hook length partition $h(\la)$ is in $D\cap O(n)$,
and then the two doubling classes of $\tilde S_n$
projecting to cycle type~$h(\la)$
each split a second time, into two classes of $\tilde A_n$.
As pointed out above,
the spin characters $\gen{\gen{h(\la)}}_+,\gen{\gen{h(\la)}}_-$
differ on the two classes  of $\tilde A_n$ contained in one
$\tilde S_n$-class but projecting
onto different $A_n$-classes of cycle type $h(\lambda)$.
Note that also the non-faithful characters $\gena{\la}_\pm$ differ on
these classes.

\medskip

The alternating groups $A_n$ and the double covers $\tilde S_n$ and
$\tilde A_n$ exhibit a special phenomenon that we want to highlight here.

\begin{df}
Let $G$ be a finite group,  $I\subset\Irr(G)$. A pair of conjugacy classes $x^G$ and $y^G$
is said to  be {\em critical} for $I$ if
$\chi(x)\neq \chi(y)$, for $\chi\in I$, but $\chi(x)=\chi(y)$ for all $\chi\in \Irr(G)\setminus I$.
\\
We call a critical pair $x^G$ and $y^G$ for $I=\{\chi_1,\chi_2\}\subset \Irr(G)$
a pair of {\em detecting classes} for $\chi_1,\chi_2$ if
$\chi_1(x)-\chi_1(y)=\chi_2(y)-\chi_2(x)$.
\end{df}

Of course, the set $I$ above should be taken of small size to give interesting information.

\begin{remexs}
{\rm
(i)
If there is a pair of detecting classes for $\chi_1,\chi_2$, then
$\chi_2=\overline{\chi_1}$, or the characters are both real.

(ii)
For $A_n$, for each symmetric partition $\mu$ of $n$,
the pair of classes $\sigma^\pm_{h(\mu)}$ is a detecting pair for $\{\mu\}_\pm$.

(iii)
For $\tilde{S}_n$, for each $\la\in D^-(n)$, the pair of classes of elements  projecting to cycle type $\la$ is a detecting pair for
$\gen{\la}_\pm$.

(iv)
For $\tilde{A}_n$, for each $\la\in D^+(n)$, the two pairs of classes of elements of $\tilde{A}_n$ projecting to cycle type $\la$ and belonging to one $\tilde{S}_n$ class are detecting pairs
for~$\gen{\gen{\la}}_\pm$.

(v)
While the situations above are the ones used in the later sections,
it should be noted that there many more such instances,
and we mention just a few examples.
\\
For $G=\GL(3,2)$, the pair of classes of elements of order~7 is
detecting for the pair of characters of degree~3.
For $G=\PSL(2,11)$, the pair of classes of elements of order~5
(order 11, resp.) is
detecting for the pair of characters of degree~12
(degree~5, resp.).
For $G=M_{11}$, the pair of classes of elements of order~8
(order 11, resp.) is
detecting for the conjugate pair of characters of degree~10
(degree 16, resp.).
}
\end{remexs}

The reason for the notions of critical and detecting classes is the following easy but very useful result on Kronecker products.
It originates with the usage of detecting classes in \cite{BB}, its variations in \cite{PPV,PP}, and the idea to consider
detecting classes for pairs of spin characters of the double
cover groups; we will follow this up in later sections.

\begin{lm}\label{lm:criticaldetect}
Let $G$ be a finite group, $x,y\in G$.
Let $\psi$ be a character of $G$ with $\psi(x)=\psi(y)\neq 0$.
\begin{enumerate}
\item[{(1)}]
Assume that $x^G,y^G$ is a critical pair for $I\subset \Irr(G)$.
Then for any $\chi\in I$, $\psi\cdot \chi$ has a constituent in $I$.
\\
In particular, if $\psi$ is irreducible and $\chi \in I$,
then $\chi \cdot \sum_{\nu\in I} \nu$
contains $\psi$ as a constituent.
\item[{(2)}]
Assume that $x^G,y^G$ is a detecting pair for $\chi_1,\chi_2\in \Irr(G)$.
Set $m_j=\gen{\psi \cdot \chi_1,\chi_j}$, $j=1,2$.
Then
$$
\psi(x)=m_1-m_2\:.
$$
Furthermore,
$$\max(m_1,m_2) \ge |\psi(x)|>0\:.$$
\end{enumerate}
\end{lm}

\proof
The assertions in (1) are obvious
(note that $I$ is closed under complex conjugation).
For the claim in (2), we compute the difference
of the values of $\psi\cdot \chi_j$ on the two classes.
Set $t:=\chi_1(x)-\chi_1(y)=\chi_2(y)-\chi_2(x)$; note that $t\ne 0$.
First we have
$$
\psi \chi_1(x)-\psi\chi_1(y)= \psi(x)(\chi_1(x)-\chi_1(y))
=\psi(x)t.
$$
On the other hand,
$\psi\chi_1=
m_1\chi_1+m_2\chi_2+\sum_{\chi \ne \chi_1,\chi_2} m_\chi \chi$,
and hence (as the pair $x^G,y^G$ is detecting for $\chi_1,\chi_2$) we obtain
$$
\psi \chi_1(x)-\psi\chi_1(y)=
m_1(\chi_1(x)-\chi_1(y))+m_2(\chi_2(x)-\chi_2(y))
=(m_1-m_2)t.
$$
As $t\ne 0$, we deduce $\psi(x)=m_1-m_2$.
The assertion $\max(m_1,m_2) \ge |\psi(x)|$ now follows immediately.
\qed

\section{The basic spin character and hooks}

Of central interest in the representation theory
of the symmetric groups
are the Kronecker coefficients $g(\la,\mu,\nu)$
appearing as expansion coefficients in the
Kronecker products
$$
[\lambda][\mu]= \sum_{\nu} g(\lambda,\mu,\nu) [\nu] \:.
$$
Using this notation, Saxl's conjecture may be rephrased as saying
the following for the staircase partition $\rho_k=(k,k-1,\ldots,2,1)$ of $n=k(k+1)/2$:
$$
g(\rho_k,\rho_k,\la)>0 \quad \text{ for all partitions }
\la \text{ of } n\:.
$$

As we will show in the following,
products of spin characters or mixed products
of a spin character and a non-faithful character
can play an important
role towards finding constituents in the square $[\rho_k]^2$.

\medskip

For the product of any ordinary character $[\lambda]$ of $S_n$
with the basic spin characters~$\gen{n}_{(\pm)}$,
Stembridge has provided an efficient combinatorial formula in~\cite{Stembridge}.
It was already observed in \cite{Bmix} that
as an immediate consequence of this formula the following
result is obtained.
\begin{prop}\label{basic-hooks}
We have the following spin product decomposition:
$$
\gen{n}\cdot \widehat{\gen{n}}
= \sum_{j=0}^{n-1} [n-j,1^j] =: \chi_{\text{hook}} .
$$
\end{prop}

Regev \cite{R-hook} showed that the sum of all hook characters
has particular easy values;
Taylor recently gave a different proof of this fact \cite{T}.
Here, we point out that this is a direct consequence
of  the proposition above and the knowledge of the values of the
basic spin characters stated in Section~\ref{prelim}:
\begin{cor} 
Let $\sigma\in S_n$ be of cycle type $\alpha$.
Then
$$\chi_{\text{hook}}(\sigma)= \left\{ \begin{array}{ll}
2^{l(\alpha)-1} & \text{if } \alpha \in O(n)\\
0 & \text{otherwise} \end{array}\right. \:.$$
\end{cor}

\medskip

For our purpose of finding constituents
in the  square of the staircase character
the following result on spin products for $\tilde S_n$
turns out to be important;
it was obtained in the context of classifying the
homogeneous spin products \cite{BKspin}.
Here, it provides the crucial link to the Saxl conjecture.

\begin{prop}\label{spinhom}
Let
$n= k(k+1)/2$,
$\rho_k=(k,k-1,\ldots,2,1)$.
Then
$$
\begin{array}{l}\gen{n}_{(\pm)} \cdot\gen{\rho_k}_{(\pm)}
=2^{a(k)}[\rho_k]
 \end{array}$$
where $a(k)$ is given by
$$a(k)=\left\{\begin{array}{ll}
\frac{k-2}2 & \mbox{if $k$ is even}\\
\frac{k-1}2 & \mbox{if } k \equiv 1 \mod 4\\
\frac{k-3}2 & \mbox{if } k \equiv 3 \mod 4\\
\end{array}
\right.
\quad .
$$
\end{prop}

\medskip

The proposition above is the key for finding new constituents
in $[\rho_k]^2 $; this is due to the following observation.

\begin{lm}\label{lm:link}
Let $[\la]$ be a constituent of
$\gen{\rho_k}_{(\pm)}
\cdot\gen{\rho_k}_{(\pm)}$, for a suitable choice of associates.
Then all constituents of \: $\chi_{\text{hook}} \cdot [\la]$ \:
are constituents of  $[\rho_k]^2$.
In particular, $[\la]$  is a constituent of  $[\rho_k]^2$.
\end{lm}

\proof
For any choice of associates we have
\begin{equation}\label{stairsquare}
\gen{n}_{(\pm)} \gen{n}_{(\pm)}
\cdot\gen{\rho_k}_{(\pm)}
\cdot\gen{\rho_k}_{(\pm)}
=2^{2a(k)} [\rho_k]^2 \:.
\end{equation}
Now the assertion follows immediately from Proposition~\ref{basic-hooks}
\qed

As a contribution towards the Saxl conjecture
we immediately obtain the following result, which was proved by
very different methods in  \cite{I} and
a weaker asymptotic version in \cite{PPV}:

\begin{cor}\label{cor:hookSaxl}
All hook characters $[n-j,1^j]$ are constituents in
$[\rho_k]^2 $.
\end{cor}

\proof
For any $\mu\in D(n)$,
$\overline{\gen{\mu}_{(\pm)}}$
is one of the characters
$\gen{\mu}_{(\pm)}$, so
the product $\gen{\mu}_{(\pm)}
\cdot\gen{\mu}_{(\pm)}$
contains $[n]$ for a suitable choice of associates.
Hence the assertion follows immediately from Lemma~\ref{lm:link}.
\qed

\begin{rem}
{\rm
Applying Lemma~\ref{lm:link} with other
constituents $[\la]$
than the trivial character,
gives further near-hook constituents;
we will see more such constituents in the next section.
For an example at this point,
using the constituents $[n-3,3]$ or $[n-3,1^3]$
of $\gen{\rho_k}_{(\pm)}\cdot {\gen{\rho_k}_{(\pm)}}$,
for a suitable choice of associates,
obtained in \cite[Theorem 3.5]{Bmix}, one finds
as constituents in $[\rho_k]^2 $ characters
to double-hooks with the smaller principal hook
being of size at most~3.
}
\end{rem}

\section{Constituents in spin products and the Saxl conjecture}

As mentioned earlier, the common theme in the character theories of the alternating
groups and the double covers of the symmetric and alternating
groups that is crucial in the applications here is the existence of
critical or detecting classes.

For example,
by using the characters of the alternating groups, it was shown in~\cite{BB}
that for a symmetric partition $\lambda$, the character $[\lambda]$
is always a constituent in its own square $[\lambda]^2$.
The idea was to use the pair of conjugacy classes in the alternating groups
that detect the characters $\{\lambda\}_{\pm}$, namely the ones
of cycle type~$h(\lambda)$. 
This was taken up by Pak, Panova and Vallejo in \cite{PPV} to provide a criterion for constituents
in $[\lambda]^2$; we give a very short argument here.

\begin{lm}\cite[``Main Lemma'']{PPV}\label{MainLemma}
Let $\lambda$ be a symmetric partition.
Let $\mu$ be a partition with $[\mu](h(\la))\ne 0$.
Then $[\mu]$ is a constituent of $[\lambda]^2$.
\end{lm}

\proof
The restriction $\chi=[\mu]\downarrow_{A_n}$ takes the same
(non-zero) values on the two classes of cycle type $h(\lambda)$,
which are critical for the pair of characters $\{\lambda\}_\pm$.
Hence by Lemma~\ref{lm:criticaldetect}
the character $\chi \cdot \{\lambda\}_+$
contains one of $\{\lambda\}_\pm$, and thus the covering product
$[\mu]\cdot[\lambda]$ contains $[\lambda]$, i.e.,
$g(\lambda,\lambda,\mu)>0$, as required.
\qed

Towards the Saxl conjecture this immediately gives:
\begin{cor}\cite{PPV}
Let $\mu$ be a partition such that $[\mu](h({\rho_k}))\ne 0$.
Then $[\mu]$ is a constituent of $[\rho_k]^2$.
\end{cor}

Unfortunately, many irreducible characters vanish on the class
of cycle type $h(\rho_k)$; this was already analysed in \cite{PPV}.
It turns out that the idea of critical detecting classes
is more powerful in the context of spin characters.

\begin{theorem}\label{spinMain}
Let $\lambda \in D(n)$, $n\ge 4$.
Let $\mu$ be a partition such that $[\mu](\la)\ne 0$.
Then $[\mu]$ is a constituent of $\gen{\lambda}_{(\pm)}\widehat{\gen{\lambda}}$, i.e.,
of $\gen{\lambda}_{(\pm)}\cdot \gen{\lambda}_{(\pm)}$
for a suitable choice of associates.
\end{theorem}

\proof
For $\lambda \in D^-(n)$, this follows immediately from
Lemma~\ref{lm:criticaldetect}(1), as  the two $\tilde{S}_n$-classes projecting
to cycle type $\lambda$ are a critical pair for the pair
of characters~$\gen{\lambda}_\pm$.

For $\lambda \in D^+(n)$,
we use spin characters of the
double covers of the alternating groups.
First assume $\lambda \in D^+(n)\setminus(D\cap O(n))$.
Then the two classes of $\tilde A_n$ projecting to cycle type $\la$
are critical for the two spin characters of $\tilde A_n$ labelled by $\la$.
Hence one of~$\gen{\gen{\la}}_\pm$ is a constituent
of 
$[\mu]\downarrow_{\tilde{A}_n} \cdot \gen{\gen{\la}}_\pm$;
thus $\gen{\lambda}$ is a constituent of
$[\mu] \cdot \gen{\la}$,
and the claim follows.
Now take $\lambda \in D\cap O(n)$; note that we then have four $\tilde{A}_n$-classes projecting to type~$\la$.
We consider non-conjugate elements
$\sigma_1,\sigma_2\in \tilde{A}_n$
that belong to the same $\tilde{S}_n$-class, projecting to
cycle type $\la$.
Recall that $\gen{\gen{\la}}_\pm$ are the only irreducible
spin characters that differ on $\sigma_1,\sigma_2$.
Now $[\mu]\downarrow_{\tilde{A}_n} \cdot \gen{\gen{\la}}_\pm$
is a spin character with different values on  $\sigma_1,\sigma_2$,
hence it must have one of $\gen{\gen{\la}}_\pm$ as a constituent,
implying the claim as before.
\qed

This leads to a powerful criterion in the context of the Saxl conjecture:
\begin{cor}\label{cor:newcriterionSaxl}
Let $\mu$ be a partition such that $[\mu](\rho_k)\ne 0$.
Then $[\mu]$ is a constituent of $[\rho_k]^2$.
\end{cor}

\proof
The assertion follows by Lemma~\ref{lm:link} and Theorem~\ref{spinMain}.
\qed

\begin{rem}{\rm
Using non-vanishing on the class $\rho_k$
produces many more constituents
than non-vanishing on the class $h({\rho_k})$!
Fortunately, the two tests may be combined to give
even more constituents.
For an illustration of this,
here are the numerical values up to $k=11$:
\begin{center}
\begin{tabular}{|c|c||c|c|c|c|c|}
\hline
\hline
&&&&& non-zero & \\
$k$ & $n$ & $p(n)$ & non-zero  & non-zero & on $\rho_k$ or  & \% \\
 & & & on $h(\rho_k)$ & on $\rho_k$ &  on $h(\rho_k)$ & \\
\hline
2 & 3 & 3 & 3 & 2 & 3 & 100  \\

3 & 6 & 11 & 5 & 6 & 9 & 81.8  \\

4 & 10 & 42 & 21 & 24 & 33 & 78.6 \\

5 & 15 & 176  & 45 & 114 & 131 & 74.4 \\

6 & 21 &792 & 231 & 524 & 607 & 76.6 \\

7 & 28 & 3718 & 573 & 2408 & 2623 & 70.5 \\

8 & 36 &17977 & 3321 & 12734 & 13567 & 75.5 \\

9 & 45 & 89134 & 9321 & 67462 & 69692 & 78.2 \\

10 & 55 & 451276 & 59091 & 370590 & 381375 & 84.5 \\

11 & 66 & 2323520 & 183989 & 2036486 & 2060003 & 88.7\\
\hline \hline
\end{tabular}
\end{center}
The non-vanishing on $\rho_k$ also seems to be a good criterion  
when compared with
Ikenmeyer's criterion which gives the constituents 
to partitions comparable to $\rho_k$ in dominance order~\cite{I}; 
in the region above, the corresponding percentage is decreasing, 
and already below 50\% at $k=9$.   
}\end{rem}

\bigskip

We want to illustrate the usefulness of the criterion given above
by showing that
in particular all characters to 2-part partitions are constituents of $[\rho_k]^2$.
Note that by the criterion which uses the value on the class to $h(\rho_k)$,
asymptotically, these constituents are found in \cite[Corollary 6.4]{PPV}; on the other hand, the constituents to 2-part partitions are also obtained by Ikenmeyer's result \cite{I}.
\\
In fact, as we will see below, by the approach via spin characters
many further constituents of $[\rho_k]^2$ are obtained,
namely the constituents in products of hook characters and characters to
2-part partitions.

\medskip

We consider a character to a 2-part partition $[n-j,j]$,
with $1\le j \le n/2$.
By the Littlewood-Richardson rule, we have
$$[n-j,j] = [n-j]\circ [j] - [n-j+1]\circ [j-1],$$
where on the right hand side the outer product of two characters
is denoted by~$\circ$.

For $k\le n$, we set
$$d_k(n)=|\{\lambda = (\ell_1 , \dots ) \in D(n) \mid \ell_1 \le k\}|\:.$$
First we want to apply the criterion given in Theorem~\ref{spinMain},
that is, we have to evaluate $[n-j,j]$ at $\rho_k$,
for $n=k(k+1)/2$. By the above, we obtain the value
$$[n-j,j](\rho_k) = d_k(j)-d_k(j-1)\:.$$
The partition numbers $d_k(m)$, $0\le m \le n=k(k+1)/2$,
are easily seen to form a symmetric sequence.
They are also the coefficients in the expansion
$$\prod_{i=1}^k (1+x^i) = \sum_{m=0}^n d_k(m) x^m\:.$$
This polynomial is known to be unimodal, by quite different and
intricate proofs due to Hughes~\cite{Hughes}, and
Odlyzko and Richmond \cite{OR};
see also Stanley's article \cite{S-unimodal} for more
on this and related unimodal sequences.
Here, based on a result by Odlyzko and Richmond \cite{OR},
we find that the sequence is in fact {\em almost strictly} unimodal;
the few exceptions for small $k$, and at the very beginning
and end (and an instance in the middle when the sequence has even length)
can be described explicitly.

\begin{theorem}\label{thm:strictunimod}
For $k\in \N$,
the following holds for the partition numbers $d_k(m)$:
\begin{enumerate}
\item[{(i)}] For $k\ge 1$, $d_k(0)=d_k(1)=1$.

\item[{(ii)}] For $k\ge 2$, $d_k(2)=1$; $d_2(3)=1$, $d_3(4)=1$.

\item[{(iii)}] For $k\ge 4$, $d_k(3)=d_k(4)=2$.

\item[{(iv)}] For $k\in \{4,\ldots,11\}$,
further equalities $d_k(m-1)=d_k(m)$ hold
for the following
values of $m \le k(k+1)/4$:
\[
\begin{array}{|c||c|c|c|c|c|c|c|c|}
\hline\hline
&&&&&&&& \\[-7pt]
k& 4&5&6&7&8&9&10&11\\
\hline
&&&&&&&& \\[-5pt]
m&5&6,7&7,8,10&8,11,13,14&16,17&19,22&26&32\\[-10pt]
&&&&&&&& \\
\hline \hline
\end{array}
\]

\item[{(v)}]
For $k\ge 12$, $d_k(m-1)<d_k(m)$
for $5\le m\le k(k+1)/4$.
\end{enumerate}

Up to equalities deduced from the above by the symmetry of the sequence (including an equality in the middle when $k \equiv 1$ or~$2\mod 4$),
all equalities in the partition sequences
$(d_k(m))_{0\le m \le  k(k+1)/2}$, $k\in \N$, are described above.
\end{theorem}

\proof
The equalities given in (i)-(iii) above are easily seen, and the exceptional
ones in (iv) are obtained by computation.
Up to the relations coming from the symmetry of the sequence, it only
remains to discuss (v) to cover all relations between $d_k(m)$ and $d_k(m-1)$.

Considering strict partitions of $m$ with largest part being smaller than $k$ or equal to $k$, respectively, gives the recursion
$$d_k(m)=d_{k-1}(m) + d_{k-1}(m-k)\:.$$
This will now be applied to show in an induction on $k$,
that for $k\ge 12$ the sequence $d_k(m)$ is strictly increasing in the range $4\le m\le k(k+1)/4$.

For $12\le k< 60$, this claim holds by direct computation.
So we assume now that $k\ge 60$.
Then, using unimodality of the sequences
and induction, for $m\le k(k-1)/4$ we obtain
$$
d_k(m)-d_k(m-1)=d_{k-1}(m)-d_{k-1}(m-1) + d_{k-1}(m-k)-d_{k-1}(m-1-k)
>0\:.
$$
For the range $k(k-1)/4 \le m \le k(k+1)/4$, the inequality
$d_k(m-1)<d_k(m)$ is provided by \cite[Theorem 4]{OR}.
Hence, from $m=4$ onwards, the $d_k$-sequence increases up to the middle;
for $k \equiv 1,2 \mod 4$, there is an equality in the middle by symmetry.
\qed

From the discussion before the theorem we now deduce immediately:
\begin{cor}
Let $k\in \N$, $n=k(k+1)/2$.
Then $[n-j,j](\rho_k)>0$ for $0\le j \le n/2$,
unless $j$ is one of the exceptional values
given in Theorem~\ref{thm:strictunimod} where $d_k(j-1)=d_k(j)$.
\end{cor}

Our criterion in Theorem~\ref{spinMain} now implies:
\begin{cor}\label{cor:2part-in-spinsquare}
Let $k\in \N$, $n=k(k+1)/2$, $0\le j \le n/2$ such that $j$
is not one of the exceptional values in Theorem~\ref{thm:strictunimod}.
Then $[n-j,j]$ is a constituent of $\gen{\rho_k}_{(\pm)}\cdot \gen{\rho_k}_{(\pm)}$
for a suitable choice of associates.
\end{cor}

Towards the Saxl conjecture we deduce:
\begin{cor}\label{cor:hook-2part-Saxl}
Let $k\in \N$, $n=k(k+1)/2$.
\begin{enumerate}
\item[{(1)}]
Assume $0\le j \le n/2$ and
$j$ is not one of the exceptional values
in Theorem~\ref{thm:strictunimod}.
Then all constituents of
$\chi_{\text{hook}} \cdot [n-j,j]$
are constituents of $[\rho_k]^2$.

\item[{(2)}]
All characters $[n-j,j]$, $0\le j \le n/2$,
are constituents  of $[\rho_k]^2$.
\end{enumerate}
\end{cor}

\proof
The first assertion follows immediately using
Lemma~\ref{lm:link}.
For the second assertion, we only have to check the
characters $[n-j,j]$ for the few exceptional values of~$j$
where  $[n-j,j](\rho_k)=0$.
For the cases $j=1, k\ge 2$, or $j=2,k\ge 3$, or $j=4,k\ge 4$,
the assertion $g(\rho_k,\rho_k,(n-j,j))>0$ is already known
by work of Saxl \cite{Saxl87}.
So it only remains to check the few exceptional
cases for $k\le 11$, which is easily done by computation.
\qed

\begin{rem}{\rm
Products of hook characters and characters to 2-part partitions
have been studied by Remmel \cite{Remmel} and Rosas \cite{Rosas}.
The extra constituents that we obtain from Corollary~\ref{cor:hook-2part-Saxl}
are all labelled by double hooks; in fact, every character to a 
double hook appears in a product of the form 
$\chi_{\text{hook}}\cdot [n-j,j]$, where there is some 
flexibility in the choice of~$j$  \cite[Theorem 2.2]{Remmel}.
Avoiding the few exceptional values for $j$ in Theorem~\ref{thm:strictunimod}
and showing $g(\rho_k,\rho_k,\la)>0$ for small $k$ by computation 
(or one of the available criteria)   
then gives all double hook characters $[\la]$ in $[\rho_k]^2$.
}\end{rem}

Already in \cite{BB}, for  symmetric $\la$
the result on the multiplicity $g(\la,\la,\la)$ was made more precise;
from the calculation of the scalar products on the level of the alternating groups a congruence mod 4 was found for the Kronecker coefficient.
Indeed, a similar calculation was used  to find the value
$|[\mu](h(\la))|$ as a lower bound
for the Kronecker coefficient $g(\la,\la,\mu)$ in \cite{PP};
we have seen a general version of this in Lemma~\ref{lm:criticaldetect}.

Also in the spin case we can make the result in Theorem~\ref{spinMain}
more precise.

\begin{theorem}\label{spinMainbound}
Let $\lambda \in D(n)$, $n\ge 4$.
Let $\mu$ be a partition such that $[\mu](\la)\ne 0$.
\begin{enumerate}
\item[{(1)}]
Let $\la \in D^-(n)$. Set\;
$m_\pm=\langle \gen{\lambda}_{\pm},[\mu]\cdot \gen{\lambda}_+\rangle$.
Then
$$m_+ - m_- = [\mu](\la)\:.$$
In particular, \;
$
0 < |[\mu](\la)| \le
\max(m_+,m_-) \le \langle [\mu],\gen{\lambda}_{\pm} \cdot \widehat{\gen{\lambda}}\rangle
\:.$
\item[{(2)}]
Let $\la \in D^+(n)$. Set\;
$m_\pm=\langle \gen{\gen{\lambda}}_{\pm},[\mu]\downarrow_{\tilde{A}_n}\cdot \gen{\gen{\lambda}}_+\rangle$.
Then
$$m_+ - m_- = [\mu](\la)$$
and
$$m_+ + m_- = 2m_- +[\mu](\la)=\gen{[\mu]\gen{\la},\gen{\la}}\:.$$
In particular,
$$
0 < |[\mu](\la)| \le \max(m_+,m_-)
\le \langle{[\mu]\gen{\la},\gen{\la}}\rangle$$
and
$$\langle{[\mu]\gen{\la},\gen{\la}}\rangle \equiv [\mu](\la) \mod 2\:.$$
\end{enumerate}
\end{theorem}

\proof
The $D^-$ case is an immediate consequence of Lemma~\ref{lm:criticaldetect}.
\\
For the $D^+$ case, we compute the difference of the values
of $[\mu]\downarrow_{\tilde{A}_n}\cdot \gen{\gen{\lambda}}_+$
on the two $\tilde{A}_n$-classes contained in one $\tilde{S}_n$ class
projecting to cycle type $\la$, similarly
as in the proof of Lemma~\ref{lm:criticaldetect}.
Note that this character is a linear combination
of irreducible spin characters, and restricted to
the spin characters,
the two $\tilde{A}_n$-classes are critical.
This gives the first assertion.
Observing that
$\langle \gen{\gen{\lambda}}_{\pm},[\mu]\downarrow_{\tilde{A}_n}\cdot \gen{\gen{\lambda}}_-\rangle = m_\mp$
then gives the second assertion.

In both cases the additional claims are an immediate consequence.
\qed

\begin{rem}
{\rm
There is a different way to obtain a character $[\mu]$
as a constituent in $[\rho_k]^2$
by using its character value on the class
of cycle type $\rho_k$.
As the criterion to be described now is much weaker, we only
provide the main arguments without going  into the
necessary background in detail; again, it is based on using a
special detection property of the class of cycle type~$\rho_k$.

By the results in \cite{BO00} and \cite{BO05},
each partition $\alpha\in O(n)$ is special
for the spin character(s) $\gen{\beta}_{(\pm)}$ to its Glaisher correspondent
$\beta\in D(n)$, with respect to the condition
that the 2-power in the
spin character value of $\gen{\beta}_{(\pm)}$ on elements
projecting to cycle type $\alpha$ is
the smallest (among the 2-powers in spin character values
on this class).
The spin character(s) to $\beta=\rho_k$ are indeed
unique (up to associates) with this property on their
special class, to the Glaisher correspondent $\alpha$ of $\rho_k$.
Thus, when we multiply $\gen{\rho_k}_{(\pm)}$
with a character $[\mu]$ that has odd value on $\alpha$,
the product has to contain $\gen{\rho_k}_{(\pm)}$
as a constituent.
Now it follows from a general character-theoretic fact \cite[(6.4)]{Feit}
that the values  $[\mu](\rho_k)$ and $[\mu](\alpha)$ are congruent
modulo~2.
Hence, whenever $[\mu](\rho_k)$ is odd,
we deduce that $[\mu]$ is a constituent
of $\gen{\rho_k}_{(\pm)}\widehat{\gen{\rho_k}}$,
and hence also of $[\rho_k]^2$.
}
\end{rem}

\section{A spin variant of Saxl's conjecture}

The product of two spin characters for the double cover groups
decomposes into
non-faithful irreducible characters,
while a mixed product of a spin character and a non-faithful
character decomposes into irreducible spin characters.
Thus in a variant of Saxl's conjecture for the double cover groups
we cannot expect to find an irreducible
character whose square contains all irreducible characters
but have to be more modest.

We start with a result obtained in the context of
classifying homogeneous mixed products; it is a special product appearing in
\cite[Theorem 3.2]{Bmix}:

\begin{prop}\label{mixhom} 
For $n=k^2\geq 4$ we have
$$ \gen{n}_{(\pm )} \cdot [k^k] =
2^{\lfloor \frac{k-1}2 \rfloor} \gen{2k-1,2k-3,\ldots,3,1}\:.$$
\end{prop}
\medskip

We denote by $\tau_k=(2k-1,2k-3,\ldots,3,1) \in D(k^2)$ the
``spin staircase'' of length~$k$; note that $\tau_k=h((k^k))$.

Arguing similarly as for Lemma~\ref{lm:link}
we immediately deduce from
Proposition~\ref{mixhom}:
\begin{lm}\label{lm:spinstair}
Let $[\la]$ be a constituent of $[k^k]^2$.
Then all constituents of \: $\chi_{\text{hook}} \cdot [\la]$ \:
are constituents of  $\gen{\tau_k}^2$.

In particular,  the square
$\gen{\tau_k}^2$ contains
$[\la]$  and all hook characters
$[n-j,1^j]$, $j\in \{0,1,\ldots,n-1\}$.
\end{lm}

Using Lemma~\ref{MainLemma},
we deduce from Lemma~\ref{lm:spinstair}:
\begin{cor}\label{cor:crit-spinSaxl}
Let $n=k^2$. Then the square
$\gen{\tau_k}^2$ contains all characters
$[\mu]$
such that $[\mu](\tau_k)\ne 0$.
\end{cor}

\begin{rem}{\rm
In fact, for $2\le k \le 5$, $\gen{\tau_k}^2$ contains
all characters~$[\mu]$, $\mu\in P(k^2)$.
So as a spin variant of Saxl's
conjecture we may ask whether this holds for all $k\ge 2$.

From the criterion in Corollary~\ref{cor:crit-spinSaxl}
above, already many constituents
are obtained; here are the numerical values up to $k=8$:
\begin{center}
\begin{tabular}{|c|c|c|c|c|}
\hline
\hline
&&&&\\[-10pt]
$k$ & $n$ & $p(n)$ & non-zero &   \% \\
 & &  & on $\tau_k$ & \\
\hline
&&&&\\[-10pt]
1 & 1 & 1 & 1 & 100  \\

2 & 4 & 5 & 3 & 60.0  \\

3 & 9 & 30 & 15 & 50.0  \\

4 & 16 & 231 & 93 & 40.3 \\

5 & 25 & 1958  & 755 & 38.6 \\

6 & 36 & 17977 & 7185 & 40.0\\

7 & 49 & 173525 & 75430 & 43.5 \\

8 & 64 &  851522 & 1741630  & 48.9 \\

\hline
\end{tabular}
\end{center}

}\end{rem}

\bigskip

Computations with GAP \cite{GAP4} led to the following
conjecture, adding to the conjectures made by Heide, Saxl,
Tiep and Zalesskii~\cite{HSTZ}
on character squares;
note that for $\la\in D^-(n)$, always one of
$[n]$ or $[1^n]$ is missing from the square $\gen{\la}^2$.

\begin{conj}
For any $n\geq 4$, $n\neq 5$, there is a spin character
$\gen{\la}$, $\la \in D^+(n)$,
whose square $\gen{\la}^2$ contains all
$[\mu]$, $\mu \in P(n)$.
\end{conj}

It was also conjectured in \cite{HSTZ}  that for $n\ge 4$
there is always
an irreducible character of the alternating group $A_n$
whose square contains all irreducible $A_n$ characters;
in fact, a quick check with GAP up to $n=25$
shows that for growing $n$ a large percentage of
irreducible $A_n$ characters has this property.
For the double cover groups~$\tilde{A}_n$, data
computed by GAP up to $n=25$ provides evidence for
the following conjecture;
again, for growing $n$, it seems that a large percentage of
irreducible $\tilde{A}_n$ spin characters
has the property considered here.

\begin{conj}
For any $n\geq 5$ there is a spin character
$\chi \in \Irr(\tilde{A}_n)$
whose square $\chi^2$ contains all
non-faithful $\psi \in \Irr(\tilde{A}_n)$.
\end{conj}


\end{document}